\documentclass[letterpaper, conference, twocolumn, final]{ieeeconf}
\IEEEoverridecommandlockouts 

\usepackage{soul}
\usepackage{mathtools}
\usepackage{bm}
\usepackage{algorithmicx}
\usepackage{algorithm}
\usepackage{algpseudocode}
\usepackage{graphicx}
\usepackage{epstopdf}
\usepackage{amsfonts}
\usepackage{amsmath}

\usepackage{mathtools}

\usepackage{xcolor}
\usepackage{amssymb}
\usepackage{pgfplots,tikz}
\usetikzlibrary{spy}

\pgfplotsset{width=10cm,compat=newest}
\usepgfplotslibrary{units}
\usetikzlibrary{backgrounds}
\usepackage{pgfplotstable}

\pgfplotstableread{
0.0    1.0
0.0   -0.5
}\datatable

\pgfplotsset{compat=1.3}
\newtheorem{theorem}{Theorem}[section]

\newtheorem{lemma}[theorem]{Lemma}

\newtheorem{assumption}[theorem]{Assumption}

\newtheorem{definition}[theorem]{Definition}

\def\QEDclosed{\mbox{\rule[0pt]{1.3ex}{1.3ex}}} 

\def\QED{\QEDclosed} 

\newcommand\oprocendsymbol{\hbox{$\square$}}
\newcommand\oprocend{\relax\ifmmode\else\unskip\hfill\fi\oprocendsymbol}

\usepackage{tikz}
\usetikzlibrary{arrows}
\usetikzlibrary{spy}
\usetikzlibrary{calc}

\tikzset{dashdot/.style={dash pattern=on .4pt off 3pt on 4pt off 3pt}}

\usetikzlibrary{plotmarks}
\usepackage{pgfplots}
\usetikzlibrary{intersections}
\usetikzlibrary{patterns}
\pgfplotsset{compat=1.11} 

\usepackage[printonlyused]{acronym}

\acrodef{MILP}{\emph{Mixed Integer Linear Programming}} 
\acrodef{MIG}{\emph{Mixed Integer Gomory}}
\acrodef{DiMILP}{\emph{Distributed MILP}} 
\newcommand{\DiMILP}{\texttt{DiMILP}}

\usepackage{soul}

\newcommand{\xsf}{u}
\newcommand{\var}{z}
\newcommand{\varint}{x}

\newcommand{\cost}{J}
\newcommand{\diam}{d_{\mathcal{G}}}

\newcommand{\real}{\mathbb{R}}
\newcommand{\integer}{\mathbb{Z}}

\newcommand{\convS}{\text{conv}(S)}

\DeclarePairedDelimiter\floor{\lfloor}{\rfloor}

\newcommand{\subj}{\text{subj. to}}

\DeclareMathOperator*{\argmin}{arg\,min}

\makeatletter
\newcommand{\StatexIndent}[1][3]{%
  \setlength\@tempdima{\algorithmicindent}%
  \Statex\hskip\dimexpr#1\@tempdima\relax}
\makeatother

\begin{document}
\title{
A Finite-Time Cutting Plane Algorithm for \\
Distributed Mixed Integer Linear Programming 
  }

\author{Andrea Testa, Alessandro Rucco, Giuseppe Notarstefano
\thanks{Andrea Testa, Alessandro Rucco, Giuseppe Notarstefano are with 
        Department of Engineering, Universit\`a del Salento, Lecce, Italy
        {\tt\small {\{name.lastname\}@unisalento.it} }}%
\thanks{
This result is part of a project that has received funding from
the European Research Council (ERC) under the European Union's
Horizon 2020 research and innovation programme (grant agreement
No 638992 - OPT4SMART).
}
} 

\maketitle

\begin{abstract}
  Many problems of interest for cyber-physical network systems can be formulated
  as Mixed Integer Linear Programs in which the constraints are distributed
  among the agents.  In this paper we propose a distributed algorithm to solve
  this class of optimization problems in a peer-to-peer network with no
  coordinator and with limited computation and communication capabilities. In
  the proposed algorithm, at each communication round, agents solve locally a
  small LP, generate suitable cutting planes, namely intersection cuts and
  cost-based cuts, and communicate a fixed number of active constraints, i.e., a
  candidate optimal basis. We prove that, if the cost is integer, the algorithm
  converges to the lexicographically minimal optimal solution in a finite number
  of communication rounds. Finally, through numerical computations, we analyze
  the algorithm convergence as a function of the network size. %
\end{abstract}
\IEEEpeerreviewmaketitle
\section{Introduction}
\ac{MILP} plays an important role in many problems in control, including
control of hybrid systems~\cite{bemporad1999control}, trajectory planning~\cite{richards2002spacecraft}, and task
assignment~\cite{bellingham2003multi}.
For example, thanks to the simultaneous presence of equality/inequality
constraints depending on some integer variables, nonlinear optimal control
problems can be approximated by means of \ac{MILP}.
In this paper, we consider a distributed optimization setup in which the
constraints of the MILP are distributed among agents of a network, and propose a
distributed algorithm to solve it.

We organize the relevant literature to our paper in two main blocks:
centralized and parallel approaches to solve MILP problems in control
applications, and distributed algorithms solving linear programs and convex 
problems that can be seen as relaxations or special versions of mixed 
integer programs.
First, a centralized Model Predictive Control (MPC) scheme for solving constrained multivariable control problems is proposed in \cite{bemporad2000piecewise} and \cite{bemporad2002model}. %
The MPC is formulated as a multi-parametric MILP which avoids solving (expensive) MILPs on-line. 
In \cite{axehill2014parametric} the authors propose an algorithm to solve parametric mixed integer quadratic and linear programs. 
The algorithm uses a branch-and-bound procedure, where relaxations are solved in the nodes of a binary search tree. 
In \cite{borrelli2006milp} the authors show how to formulate a 
centralized trajectory optimization problem for multiple UAVs to a finite dimensional MILP
which is solved by using a commercial branch and bound algorithm.
In \cite{chopra2015spatio} the authors address the
multi-robot routing problem under connectivity constraints. 
The authors show that such a routing problem can be formulated as an integer
program with binary variables, and then solve %
its LP relaxation. 

As for parallel methods, in \cite{kim2013scalable} a Lagrange relaxation
approach is used in order to solve the overall MILP through a master-subproblem
architecture.  The proposed solution is applied to the demand response control
problem in smart grids. Although processors are spatially distributed, the
computation is parallel since it makes use of a central coordinator.
In \cite{vujanic2016decomposition} a dual decomposition technique is proposed
for the charging control problem of electric vehicles. Here an aggregator is
required in order to assign charging slots to each individual electric vehicle.

Second, regarding distributed optimization algorithms, we concentrate on schemes
solving linear programs and convex programs that represent a relaxation of 
suitable mixed-integer programs. 
In \cite{richert2016distributed_a} the authors design a robust, distributed
algorithm to solve linear programs over networks with event-triggered
communication. 
Based on state-based rules, the agents decide when to broadcast state
information to their neighbors in order to ensure asymptotic convergence to a
solution of the linear program.
In \cite{richert2016distributed_b} the authors propose a distributed algorithm
to find valid solutions for the bargaining problem by means of a LP relaxation. 
In \cite{fischione2011utility} the authors address the Utility Maximization
problem which is, in its general formulation, a mixed-integer nonlinear
programming problem. The proposed solution is based on a convex relaxation,
i.e., the integer constraint on the rates is relaxed thus yielding a convex
program.
In \cite{wei2013distributed_alg} and \cite{wei2013distributed_conv} the authors propose a Newton-type fast converging algorithm to solve, 
under the assumption that the utility functions are self-concordant, 
the Network Utility Maximization problem. %
In \cite{notarstefano2011distributed} the authors propose constraints consensus
algorithms to solve abstract optimization programs (i.e., a generalization of
linear programs).
A distributed simplex algorithm is proposed in \cite{burger2012distributed} to solve degenerate LPs and multi-agent assignment problems in asynchronous networks. 
A distributed version of the Hungarian method is proposed
in~\cite{chopra2017distributed} to solve distributed LP arising in multi-robot
assignment problems.
In \cite{franceschelli2013gossip} the authors address multi-agent task
assignment and routing problems, modeled as MILP, in a distributed fashion. A
gossip algorithm exploiting pairwise task exchanges between agents is proposed
to find a common feasible assignment.
Finally, a distributed trajectory optimization algorithm for cooperative UAVs is
proposed in \cite{kuwata2011cooperative}.  The algorithm is based on a special
sequential computation set-up in which local \ac{MILP}s are solved by the UAVs
in a given sequence. %

The main contribution of this paper is the design of a \emph{Distributed MILP}
algorithm, called \DiMILP, that, solves \ac{MILP} problems in finite time under
the assumption of integer optimal cost. To the best of our knowledge, this is
the first distributed algorithm solving MILP problems in asynchronous, directed
networks.
The algorithm is based on the local generation of suitable cutting planes,
namely intersection cuts and cost-based cuts, and the exchange of active (basic)
constraints. %
We consider a peer-to-peer network with no coordinator, in which each agent
performs simple computations, i.e., solves small LPs and generates cutting
planes, and communicates with other agents only a small number of linear
constraints, i.e., a basis of its local LP relaxation.
By exploiting the structure of intersection cuts and cost-based cuts, we prove 
the correctness of the proposed algorithm and its convergence in a finite 
number of communication rounds. 
We analyze and discuss a set of simulations to study the evolution and the
convergence of the algorithm while varying the network size. %

We highlight some meaningful differences with respect to the literature discussed above. 
Although constraint exchange and cutting plane approaches have been proposed in
\cite{notarstefano2011distributed} and \cite{burger2014polyhedral}, 
in this paper we consider a MILP optimization problem.  The
presence of variables subject to integer constraints gives raise to new
challenges in the algorithm design and analysis.
In \cite{franceschelli2013gossip} and~\cite{kuwata2011cooperative}, the solution
to each local MILP is optimal, yet there is no guarantee on the global
optimality. On the contrary, we propose an algorithm for distributed MILP with
guaranteed finite-time convergence to a global optimum.  

The paper is organized as follows. 
In Section~\ref{sec:pf}, we introduce the MILP problem and its distributed formulation. 
In Section~\ref{sec:MILP}, we describe the cutting plane approach for \ac{MILP}. 
The distributed algorithm is introduced and analyzed in Section~\ref{sec:DiMILP}. 
Numerical computations are provided in Section~\ref{sec:NumComp} followed by the conclusion in Section~\ref{sec:concl}. %

\section{Problem set-up}
\label{sec:pf}
We consider the following MILP: 
\begin{equation} \label{eq:MILP}
	\begin{split}
		\min_{\var} &\; c^\top \var\\ 
		\subj &\; a_i^\top \var \leq b_i \,, \,i = 1, \ldots, n\\
		&\; \var \in \integer^{d_Z} \times \real^{d_R} \\
	\end{split}
\end{equation}
where $a_i \in \real^{d}$, $b_i \in \real$, $c \in \real^{d}$, and 
$n$ is the number of inequality constraints. 
Before formulating the distributed optimization set-up considered in the paper, 
we provide some useful notation.

\paragraph*{Notation}
We denote by $\var$ the decision variables in
$\integer^{d_Z} \times \real^{d_R}$, by $\varint$ the variables in
$\integer^{d_Z}$, i.e., the ones subject to integer constraints, and by $y$ the
variables in $\real^{d_R}$. %
We let $d = d_Z + d_R$.
Given an inequality $a^\top \var \leq b$ for $\var\in\real^{d} $, with $a\in\real^d$ and $b\in\real$, we use the following
simplified notation $\{ a^\top z \leq b\} := \{ z \in \real^d : a^\top z \leq b\}$ for
the related half-space. 
The polyhedron\footnote{A polyhedron is a set described by the intersection of a finite number of half-spaces.} 
induced by the inequality constraints
$a_i^\top \var \leq b_i \,, \,i \in \{1, \ldots, n\}$, is
$P:=\bigcap_{i=1}^n \{ a_i \var \leq b_i \}$.
Given two vectors $v, w \in \real^d$, $v$ is lexicographically greater than $w$, $v >_{lex} w$, if there exists $k \in \{1,\ldots, d\}$ such that $v_k > w_k$ and $v_m = w_m$ for all $m<k$.  
  
In this paper we assume an LP solver is available. In particular, we consider a
solver that is able to find the unique \emph{lexicographically minimal optimal
  solution} of the problem. 
  From now on, we call such a solver \textsc{LPlexsolv} and say that
it returns the \emph{lex-optimal} solution meaning it is the lexicographically
minimal optimal solution of the solved LP problem. 
\textsc{LPlexsolv} also returns an optimal \emph{basis} identifying the lex-optimal solution.
Given an LP problem with constraint set $P:=\bigcap_{i=1}^n P_i$, with each
$P_i$ a half-space, a basis $B$ is the intersection of a \emph{minimal} number
of half-spaces $P_{\ell_1}, \ldots,P_{\ell_k}$, $k\leq d$, such that the
solution of the LP over the constraint set $B$ is the same as the one over
$P$. If the lex-optimal solution is considered, it turns out that $B$ is the
intersection of exactly $d$ half-spaces ($k=d$).

\medskip

In our distributed setup, we consider a network composed by a set of agents $V=\{1, \ldots, N_{ag}\}$. 
In general, the $n\geq N_{ag}$ constraints in Problem~\ref{eq:MILP} are
distributed among the agents, so that each agent knows only a small number of
constraints. 
For simplicity, we assume one constraint $\{a_i^\top z \leq b_i\}$ is assigned
to the $i$-th agent, so that $N_{ag} = n$, but we will keep the two notations
separate to show that the algorithm can be easily implemented also when
$N_{ag}<n$.  
The communication among the agents is modeled by a 
time-varying digraph $\mathcal{G}_c(t)=(V,E(t))$, with $t\in\mathbb{N}$ being a
universal slotted time. A digraph $\mathcal{G}_c(t)$ models the communication in
the sense that there is an edge $(i,j) \in E(t)$ if and only if
agent $i$ is able to send information to agent $j$ at time $t$.
For each node $i$, the set of \emph{in-neighbors} of $i$ at time $t$ is denoted
by $N_i(t)$ and is the set of $j$ such that there exists an edge $(j,i) \in E(t)$.
A static digraph is said to be \emph{strongly connected} if there exist a
directed path for each pair of agents $i$ and $j$.  For a time-varying
digraph, we require the \emph{joint strong connectivity}, i.e.,
$\forall t \in \mathbb{N}, \cup_{\tau=t}^{\infty} \mathcal{G}_c(\tau)$ is strongly
connected.

\section{Cutting Planes in Mixed Integer Linear Programming}
\label{sec:MILP}
In this section we provide a brief description of one of the most used methods
to solve a centralized MILP, i.e., the cutting plane approach. %
\vspace{-0.1cm}
\subsection{Cutting-Plane approach for MILP}
\begin{figure}[htpb]
   \begin{center}
\begin{tikzpicture}[
    scale=0.85,
    axis/.style={->, >=stealth'},
    important line/.style={thick},
    dashed line/.style={dashed, thin},
    pile/.style={thick, ->, >=stealth', shorten <=2pt, shorten
    >=1pt},
    every node/.style={color=black}
    ]
\begin{axis}[
    xmin=-4,xmax=4,
    ymin=-3,ymax=4,
    grid=both,
    minor tick num=1,
    xlabel={$x$},
    ylabel={$y$},
    yticklabels={},
    xticklabels={},
    axis x line=middle, axis y line=middle]
        \addplot[green,domain=-5:5,samples=10]
        { (18 - 1*x)/(10) };   
        \draw[line width=0.25mm, green]  (-2.6897, 2.069) -- node[above,sloped] {$x+10y\leq18$} (1.333, 1.6667);
        \addplot[green,domain=-5:5,samples=10]
        { (4 + 5*x)/(-8) };   
        \draw[line width=0.25mm, green]  (-2.3158, 0.9474) -- node[below,sloped] {$-5x-8y\leq4$} (0.8,-1);
        \addplot[green,domain=-5:5,samples=10]
        { (18 + 9*x)/(-3) };   
        \draw[line width=0.25mm, green]  (-2.6897, 2.069) -- node[below,sloped] {$-9x-3y\leq18$} (-2.3158, 0.9474);
        \addplot[green,domain=-5:5,samples=10]
        { (10 - 10*x)/(-2) };   
        \draw[line width=0.25mm, green]  (0.8, -1) -- node[below,sloped] {$10x-2y\leq10$} (1.333, 1.6667);
%
        \fill[green, opacity=0.5] (-2.6897, 2.069) -- (-2.3158, 0.9474)  -- (0.8, -1) --  (1.333, 1.6667) -- cycle; 
        \fill[red, opacity=0.5] (-2,2) -- (-2,0.75) -- (-1, 0.125) -- (0,-0.5) -- (1,0) -- (1,1.7) -- cycle;
        \draw[line width=0.4mm, red]  (-2, 2) -- (-2, 0.75);
        \draw[line width=0.4mm, red]  (-1, 0.125) -- (-1, 1.9);
        \draw[line width=0.4mm, red]  (0, -0.5) -- (0, 1.8);
        \draw[line width=0.4mm, red]  (1, 0) -- (1, 1.7);
    \end{axis}
\end{tikzpicture}
   \end{center}
   \caption{Mixed integer set for a MILP with $x \in \integer$, $y \in \real$. The
  polyhedron $P$ (green area) induced by the four inequality constraints (green
  lines), the set $S$ of feasible solutions of the \!MILP (red lines), the
  convex hull of $S$ (red area) are shown.  }
	\label{fig:PSconvS}
 \end{figure}
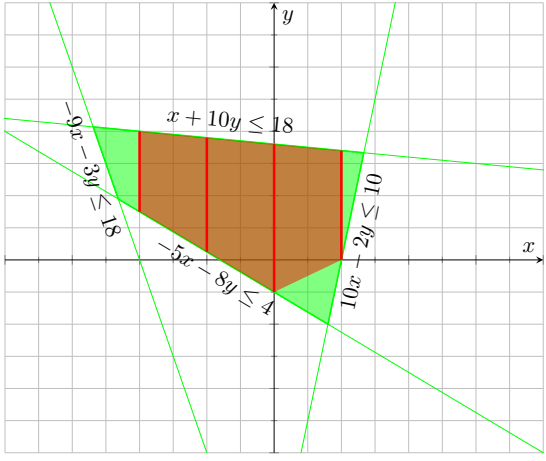

Let $S:=P \cap (\integer^{d_Z} \times \real^{d_R})$, problem \eqref{eq:MILP} is
equivalent, \cite{schrijver1998theory}, to the following LP
\begin{equation} \label{eq:MILP2}
\begin{split}
\min_{z} &\; c^\top \var \\ 
\subj &\; \var \in \text{conv}(S) \\
\end{split}
\end{equation}
where $\text{conv}(S)$ is the \emph{convex hull} of $S$. A two dimensional
example is shown in Figure~\ref{fig:PSconvS}. %
It is worth noting that, if $P$ is a bounded polyhedron, by Meyer's
Theorem,~\cite{junger200950}, $\text{conv}(S)$ is a polyhedron (i.e.,
$\text{conv}(S)$ is the solution set of a finite system of linear
inequalities).  For this reason, we make the following assumption, which is
common in MILP literature. 
\begin{assumption}[Boundedness of $P$] \label{ass:Pbound}
	The polyhedron $P$ is bounded and therefore $\convS$ is a polyhedron. 
	\oprocend
\end{assumption}

The main idea of the cutting plane approach is to neglect the integer constraints on the decision variables $\varint$, 
solve the \emph{relaxed} linear problem (i.e., with $\varint \in \real^{d_Z}$)
and properly tighten (in an iterative manner) the polyhedron $P$ until the solution is in
$\text{conv}(S)$.
The cutting plane procedure for MILPs can be summarized as follows.
Relax the integer constraint in the formulation~\eqref{eq:MILP} and solve the optimization problem by using \textsc{LPlexsolv}. 
Let $\var_{LP}$ be the lex-optimal solution. If $\var_{LP} \in \text{conv}(S)$, then $\var_{LP}$ is the lex-optimal solution of~\eqref{eq:MILP2} and, therefore, of~\eqref{eq:MILP}. 
If $\var_{LP} \notin \text{conv}(S)$, find a linear inequality (called \emph{Cutting Plane})
\begin{equation} \label{eq:CP}
	\alpha^\top \var \leq \alpha_0 \,,
\end{equation}
where $\alpha \in \real^d$ and $\alpha_0 \in \real$, which is satisfied by all $\var \in \text{conv}(S)$ and excluding $\var_{LP}$. %
Then, update $P$ with the new linear inequality~\eqref{eq:CP} and repeat this approach until the lex-optimal feasible solution for the MILP~\eqref{eq:MILP} has been found.

Regarding the convergence property of this approach, under suitable assumptions,
cutting plane algorithms obtain convergence after a finite number of iterations.
This is the case, for example, of Integer Linear Programs (ILPs) when all the
matrices are rational,~\cite{junger200950}, and of mixed binary
programs~\cite{cornuejols2008valid}.
For MILPs, if the optimal objective function value is integer, the first cutting
plane algorithm converging in a finite number of iterations has been proposed in
\cite{gomory1960algorithm}. %
These considerations justify the following assumption. 
\begin{assumption}[Feasibility and integer optimal cost]
	\label{ass:feasibility}
	Let $J(z)$ be the objective function, i.e., $J(z):=c^\top z$. There
        exists a lexicographically-minimal optimal solution $z^\star \in \convS$
        such that $J(z^\star) \leq J(z) ,\,\; \forall z \in \convS$. The optimal
        cost $J^\star:= J(z^\star)$ is integer-valued. \oprocend 
\end{assumption}

\subsection{Cutting-Plane via Intersection Cuts}
Many approaches have been developed to generate valid cutting planes, see \cite{junger200950} for a survey. 
Next we introduce the notions of split disjunction, %
and intersection cut, %
that will be used to
characterize the cuts we use in this paper, i.e., \ac{MIG} cuts.

\begin{definition}[Split Disjunction 
\cite{cornuejols2008valid}]
	Given $\pi \in \integer^{d_Z}$ and $\pi_0 \in \integer$, a split disjunction $D(\pi, \pi_0)$ is a set of the form 
	$D(\pi, \pi_0) := \{\pi^\top x\leq \pi_0\} \cup \{ \pi^\top x \geq \pi_0+1 \}$. \oprocend 
\end{definition}
Let $B_{LP}$ be a basis of the lex-optimal solution $\var_{LP}=(x_{LP},y_{LP})$,
for a given LP relaxation of \eqref{eq:MILP}, and 
$D(\pi,\pi_0)$ a disjunction
with respect to $x_{LP}$. Let $C(\var_{LP})$ be the translated (simplicial) cone
formed by the intersection of the halfplanes defining $B_{LP}$ and having apex
in $\var_{LP}$\footnote{Given a cone $S\subset\real^d$ and a point
  $p\in\real^d$, the set $p + S$ is a translated cone with apex in $p$. }.
An \emph{intersection cut}, can be derived by considering the intersections
between the extreme rays of $C(\var_{LP})$ and the hyperplanes defining
$D(\pi,\pi_0)$. 
A more detailed definition can be found in
\cite{balas2013generalized}. 

 In this paper we use Mixed-Integer Gomory (\ac{MIG}) cuts,
 \cite{gomory1960algorithm}, as cutting planes for our distributed algorithm.
 As shown in Appendix, \ac{MIG} cuts can be obtained by working on the tableau
 of a problem in the form \eqref{eq:SFMILP} equivalent to a given LP relaxation
 of \eqref{eq:MILP}.
 Specifically, the \ac{MIG} cut with respect to the $k$-th
 row of the tableau, expressed in terms of the decision variable $\var$, is
 given by:
\begin{equation} \label{eq:MILPCut_z}
\begin{split}
h_{\textsc{mig}} \!:=\!\! \Bigg\{ \!\sum_{\substack{\ell\in N_+}}{\!\!\bar{a}}_{k\ell}[b_B\!-\!A_B\var]_\ell \! - \! \bar{f_0} \!\!\sum_{\substack{\ell\in N_-}}{\!\!\bar{a}}_{k\ell}[b_B\!-\!A_B\var]_\ell \geq \! f_0 \! \Bigg\},\\
\end{split}
\end{equation}
where $A_B$ and $b_B$ define the constraints of the basis $B_{LP}$ (as in
equation~\eqref{eq:basisLP} in Appendix) and
$\bar{a}_{k\ell}=[A_B^{-1}]_{k\ell}$, $f_0=[A_B^{-1}b_B-\floor{A_B^{-1}b_B}]_k$,
$\bar{f_0}=\frac{f_0}{1-f_0}$, $N_+:=\{\ell: \bar{a}_{k\ell} \geq 0\}$ and
$N_-:=\{\ell:\bar{a}_{k\ell} < 0\}$. Here we use the notation $[\cdot]_{k\ell}$
to indicate the $(k,\ell)$-th element of a matrix and $[\cdot]_{k}$ to indicate
the $k$-th element of the vector inside the brackets. Details on how to generate
such a \ac{MIG} cut are given in Appendix.

As recalled in Theorem~\ref{thm:MIG_intcut} in Appendix, a \ac{MIG} cut with
respect to the $k$-th row (of problem~\eqref{eq:SFMILP}) is an intersection cut
to the split disjunction $D(e_{k},\floor{x_{{LP}_k}})$ and the basis $B_{LP}$,
with $e_k$ being the $k$-th vector of the canonical basis (e.g.,
$e_1 = [1\, 0\, \ldots\, 0]$) and $x_{{LP}_k}$ the $k$-th component of $x_{LP}$.
Let us consider now the first non-integer component of $\var_{LP}$, namely
$x_{LP_{k^{lex}}}$ where $k^{lex} = \argmin \{ k = 1, \ldots, d_{Z} : x_{{LP}_k} \notin
\integer \}$. 
We call \textsc{MIGoracle} the oracle that generates \ac{MIG}
cut~\eqref{eq:MILPCut_z} for $k=k^{lex}$.

In addition to the MIG cut, we also consider a constraint based
on the actual cost function value $h_c:=\{c^\top z \geq \lceil c^\top \var_{LP} \rceil \}$. 
Notice that, by Assumption~\ref{ass:feasibility}, $h_c$ does not cut off any solution of $\convS$. 
We refer to this inequality constraint as \emph{cost-based cut}.

Next, we recall a centralized algorithm based on \ac{MIG} and cost-based
cuts, which is a reformulation of Gomory's cutting plane algorithm,
\cite{gomory1960algorithm}, for MILPs in the form \eqref{eq:MILP}. This
version is presented, e.g., in~\cite{jorg2008k}. 
A pseudocode description of this algorithm is given in the table below, Algorithm~\ref{alg:centralized}. %
\begin{algorithm}%
	\caption{Cutting Plane Algorithm for MILP (\cite{gomory1960algorithm})}\label{alg:centralized}
	\begin{algorithmic}[0]	
		\State \textbf{Input} $P$, $c$
		\StatexIndent[0.25] $(z_{LP}, B_{LP}) =$\Call{LPlexsolv}{$P$, $c$}
		\StatexIndent[0.25] $h_c = \{c^\top z \geq \lceil c^\top z_{LP} \rceil \}$		
		\StatexIndent[0.25] \emph{while} $z_{LP} \not\in \convS$ \emph{do}	
		\StatexIndent[0.75] $h_{\textsc{mig}} =$ \Call{MIGoracle}{$z_{LP}$, $B_{LP}$}
		\StatexIndent[0.75] $P = P \cap h_{\textsc{mig}} \cap h_c$
		\StatexIndent[0.75] $(z_{LP}, B_{LP}) =$\Call{LPlexsolv}{$P$, $c$}
		\StatexIndent[0.75] $h_c = \{c^\top z \geq \lceil c^\top z_{LP} \rceil \}$
		\State \textbf{Output} $z^\star = z_{LP}$ 
	\end{algorithmic}
\end{algorithm}

It is worth noting that, at each iteration, Algorithm~\ref{alg:centralized} uses
the entire set of inequality constraints, $P$, and all the cuts generated up to
that iteration.

\section{Distributed MILP} %
\label{sec:DiMILP}
In this section we propose a Distributed MILP algorithm, called \DiMILP, based
on the local generation of cutting planes and the exchange of active
constraints. Then we prove its convergence in a finite number of communication
rounds under the assumption that the optimal cost is integer.

In contrast to the centralized approach, in the distributed setup, at the first
iterations, some agents may not have enough information to properly execute the
algorithm (e.g., only one constraint has been assigned to the agent and the
local MILP problem is unbounded).
For this reason, we initialize the algorithm by assigning to each agent a set of
artificial constraints which are inactive for the global
problem~\eqref{eq:MILP}. This method is often referred to as big-M method.
Specifically, the decision variable of each agent is delimited by a box
constraint.  In particular, for a given, sufficiently
large $M>0$, we define the bounding box
$H_M:= \bigcap\limits_{k=1}^{d} \left( \left\{ \var_k \leq{M} \right\} \cap
  \left\{ \var_k \geq{-M} \right\} \right) $.
\subsection{Algorithm description}
We now describe the proposed distributed algorithm. 

Besides the initial constraint
$h^{[i]} := \{ a_i \var \leq b_i \} \cap H_M$, 
each agent $i$ has two local states, namely $z^{[i]}$, associated to the
decision variable of MILP~\eqref{eq:MILP}, and $B^{[i]}$ being a candidate basis of the problem. 
We use the subscript $k$ to denote the $k$-th component of the local state, i.e., $z_k^{[i]}$. 

At the generic time $t$, the $i$-th agent solves a linear program in which the
common objective function $c^{\top}\var$ is minimized subject to the following
constraints: the intersection of its neighbors' candidate bases,
$\bigcap_{j \in N_i} B^{[j]}(t)$, its own candidate basis, $B^{[i]}(t)$, the
inequality constraint set at the initialization step, $h^{[i]}$, and the
cost-based cut $h_c$ obtained by rounding up the current cost.
Then, agent $i$ generates a \ac{MIG} cut based on the current (local)
lex-optimal solution. %
Finally, through a pivoting routine, named \textsc{Pivot}, the agent updates its
basis and thus the corresponding solution. %

Summing up, at each communication round, each agent has to generate a constraint
based on the local optimal cost value, share the basis with its neighbors
through local communication (according to the communication graph), solve the LP
problem, generate a \ac{MIG} cut, and update its state.
This procedure is formalized in Algorithm~\ref{alg:distributed3}, where we
dropped the dependence on $t$ to highlight that agent $i$ does not need to
know it to perform the update.
\begin{algorithm}[htpb]
  \caption{\DiMILP}\label{alg:distributed3}
  \begin{algorithmic}[0]		
    \State \textbf{State} 
    \StatexIndent[0.25]  $(z^{[i]}, B^{[i]})$
    \State \textbf{Initialization}
    \StatexIndent[0.25]  $h^{[i]} = \{ a_i z \leq b_i \} \cap H_M$
    \StatexIndent[0.25]  $(z^{[i]}, B^{[i]}) =$ \Call{LPlexsolv}{$h^{[i]}$, $c$}
    \State  \textbf{Evolution}
    \StatexIndent[0.25]  $h_c = \{ c^\top z \geq \lceil c^\top z^{[i]} \rceil \}$
    \StatexIndent[0.25]  $H_{tmp}\!=\! \left( \bigcap_{j \in N_i} B^{[j]} \right) \cap B^{[i]} \cap h^{[i]} \cap h_c$
    \StatexIndent[0.25]  $(z_{LP}, B_{LP}) =$ \Call{LPlexsolv}{$H_{tmp}$, $c$}
    \StatexIndent[0.25]  $h_{\textsc{mig}} =$ \Call{MIGoracle}{$z_{LP}$, $B_{LP}$}
    \StatexIndent[0.25]  $(z^{[i]}, B^{[i]}) =$ \Call{Pivot}{$B_{LP} \cap h_{\textsc{mig}}$, $c$}
  \end{algorithmic}
\end{algorithm}

We highlight that the proposed distributed algorithm is scalable in terms of
local memory, computation and communication.
Indeed, an agent sends to neighbors a candidate basis, which is a collection of
$d$ linear constraints.
Consistently, it receives a number of bases equal to its in-degree. Also, in the
computation it considers only two more inequality constraints at each iteration
(i.e., the \ac{MIG} cut, $h_{\textsc{mig}}$, and one constraint, $h_c$, based on
the local optimal cost-value).

\subsection{Algorithm Analysis}
The finite-time convergence and the correctness of the algorithm can be proven
in three steps. For the sake of space, the proofs are omitted in this paper and
will be provided in a forthcoming document.
First, we show that for each agent the local cost and the local state converge
in finite time.
Second, we prove that consensus among all the agents is attained for the cost
functions and for the candidate lex-optimal solutions.
Finally, we show that the common candidate solution is indeed the lex-optimal
solution of the global problem.

From now on, we denote $J^{[i]}(t)$ the local cost function value related to the
decision variable of agent $i$, i.e., $J^{[i]}(t) = c^\top \var^{[i]}(t)$. 

\begin{lemma}[Local convergence] \label{lemma:conv} Let $\var^{[i]}(t)$
  be the local candidate lex-optimal solution and $\cost(\var^{[i]}(t))$ the associated cost
  of agent $i$ at time $t\geq 0$ executing \DiMILP. 
  Then, in a finite number of communication rounds:
\begin{enumerate}
  \item[i)] the sequence $\{\cost(\var^{[i]}(t))\}_{t \geq 0}$ converges to a constant
  value $\bar{\cost}^{[i]}$, and
  \item[ii)] the sequence $\{z^{[i]}(t)\}_{t \geq 0}$ converges to
  $\bar{z}^{[i]} = (\bar{x}^{[i]}, \bar{y}^{[i]})$, where
  $\bar{x}^{[i]} \in \integer^{d_Z}$.
\end{enumerate} \oprocend
\end{lemma}

\begin{lemma}[Consensus] \label{lemma:agreement} Assume the communication
  network, $\mathcal{G}_c(t)$, is jointly strongly connected. Then,
  $\bar{\cost}^{[i]} = \bar{\cost}^{[j]}$ and $\bar{z}^{[i]}=\bar{z}^{[j]}$ for
  all $i,j \in \{ 1, \cdots, N_{ag}\}$. \oprocend
\end{lemma}

We are now ready to present the main result of the paper. 
\begin{theorem}[\DiMILP~convergence] %
  Consider MILP problem~\eqref{eq:MILP} in which the constraints are distributed
  among agents communicating according to a jointly strongly connected
  communication graph, $\mathcal{G}_c(t)$,
  $t\geq0$. %
  Let Assumption~\ref{ass:feasibility} hold and 
  $M>0$ be sufficiently large such that the lex-optimal solution
  of~\eqref{eq:MILP} does not change if the constraint set
  $\bigcap\limits_{i=1}^{n}h^{[i]}$ is replaced by
  $\left( \bigcap\limits_{i=1}^{n} h^{[i]} \right) \bigcap H_M$.
  Then \DiMILP\/ solves MILP problem~\eqref{eq:MILP} in a finite number of
  communication rounds. That is, the sequences
  $\{\cost(\var^{[i]}(t))\}_{t \geq 0}$ and $\{\var^{[i]}(t)\}_{t \geq 0}$,
  $i = \{ 1, \ldots N_{ag}\}$, converge, respectively, to the global optimal
  cost and to the lex-optimal solution of problem~\eqref{eq:MILP}. \oprocend
\end{theorem} 

\section{Numerical Computations}
\label{sec:NumComp}
In this section we provide numerical computations showing the effectiveness of the proposed algorithm. 

\begin{figure}[t!]
	\centering
	\input{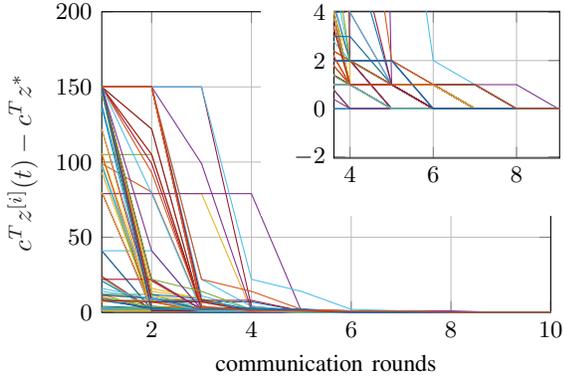}	
	\caption{Difference between the optimal cost and the cost evaluated by each agents at each communication round. 
	The agents reach the optimal cost after $9$ communication rounds. %
	}
	\label{fig:dmilp100fval}
\end{figure}

We randomly generate the MILP data as follows.  
We consider a two-dimensional problem, $d=2$. The decision variable is $z = (x, y)$ where $x \in \integer$ and $y \in \real$. 
We consider $N_{ag} = 100$ and $n = 100$, i.e., each agent only knows one constraint of the centralized MILP. 
The inequality constraints are randomly generated from the standard Gaussian
distribution (we check feasibility and discard infeasible instances). 
The cost function is $c=[1, \; 0]^\top$. 
We consider an Erd\H{o}s-R\'{e}nyi static digraph with parameter
0.015. %
We run the algorithm by setting the bounding box $H_M$ with $M = 150$. 
In Figure~\ref{fig:dmilp100fval} we show the difference between the optimal costs of each agent and global optimal cost found by solving the centralized MILP (we use the ``intlinprog'' function in
MATLAB). %
We obtain convergence to the global optimal solution and the corresponding optimal cost after $9$ communication rounds, see zoom-in in Figure~\ref{fig:dmilp100fval}. %

\begin{figure}[t!]
	\centering
%
%
\begin{tikzpicture}[font=\small]

\begin{axis}[%
width=5.987cm,
height=4cm,
at={(0cm,0cm)},
scale only axis,
xmin=4,
xmax=68,
xtick={8,16,32,64},
xlabel={number of agents},
ymin=9.6,
ymax=282.4,
ytick={  0,  25,  50,  75, 100, 125, 150, 175, 200, 225, 250, 275},
ylabel style={},
ylabel={communication rounds},
axis background/.style={fill=white},
xmajorgrids,
ymajorgrids
]
\addplot [color=black, dashed, forget plot]
  table[row sep=crcr]{%
8	29\\
8	34\\
};
\addplot [color=black, dashed, forget plot]
  table[row sep=crcr]{%
16	60\\
16	70\\
};
\addplot [color=black, dashed, forget plot]
  table[row sep=crcr]{%
32	122\\
32	143\\
};
\addplot [color=black, dashed, forget plot]
  table[row sep=crcr]{%
64	241\\
64	270\\
};
\addplot [color=black, dashed, forget plot]
  table[row sep=crcr]{%
8	22\\
8	24\\
};
\addplot [color=black, dashed, forget plot]
  table[row sep=crcr]{%
16	47\\
16	49\\
};
\addplot [color=black, dashed, forget plot]
  table[row sep=crcr]{%
32	94\\
32	96\\
};
\addplot [color=black, dashed, forget plot]
  table[row sep=crcr]{%
64	191\\
64	197\\
};
\addplot [color=black, forget plot]
  table[row sep=crcr]{%
7	34\\
9	34\\
};
\addplot [color=black, forget plot]
  table[row sep=crcr]{%
15	70\\
17	70\\
};
\addplot [color=black, forget plot]
  table[row sep=crcr]{%
31	143\\
33	143\\
};
\addplot [color=black, forget plot]
  table[row sep=crcr]{%
63	270\\
65	270\\
};
\addplot [color=black, forget plot]
  table[row sep=crcr]{%
7	22\\
9	22\\
};
\addplot [color=black, forget plot]
  table[row sep=crcr]{%
15	47\\
17	47\\
};
\addplot [color=black, forget plot]
  table[row sep=crcr]{%
31	94\\
33	94\\
};
\addplot [color=black, forget plot]
  table[row sep=crcr]{%
63	191\\
65	191\\
};
\addplot [color=blue, forget plot]
  table[row sep=crcr]{%
6	24\\
6	29\\
10	29\\
10	24\\
6	24\\
};
\addplot [color=blue, forget plot]
  table[row sep=crcr]{%
14	49\\
14	60\\
18	60\\
18	49\\
14	49\\
};
\addplot [color=blue, forget plot]
  table[row sep=crcr]{%
30	96\\
30	122\\
34	122\\
34	96\\
30	96\\
};
\addplot [color=blue, forget plot]
  table[row sep=crcr]{%
62	197\\
62	241\\
66	241\\
66	197\\
62	197\\
};
\addplot [color=red, forget plot]
  table[row sep=crcr]{%
6	26.5\\
10	26.5\\
};
\addplot [color=red, forget plot]
  table[row sep=crcr]{%
14	55\\
18	55\\
};
\addplot [color=red, forget plot]
  table[row sep=crcr]{%
30	110.5\\
34	110.5\\
};
\addplot [color=red, forget plot]
  table[row sep=crcr]{%
62	213.5\\
66	213.5\\
};
\end{axis}
\end{tikzpicture}%
	\caption{Communication rounds evolution while varying the graph diameter. 
	Box plot shows the minimum and maximum communication rounds (whiskers), 25\% and 75\% percentiles (lower and upper limit of the box) and median (red line). 
	}
	\label{fig:boxplot}
\end{figure}
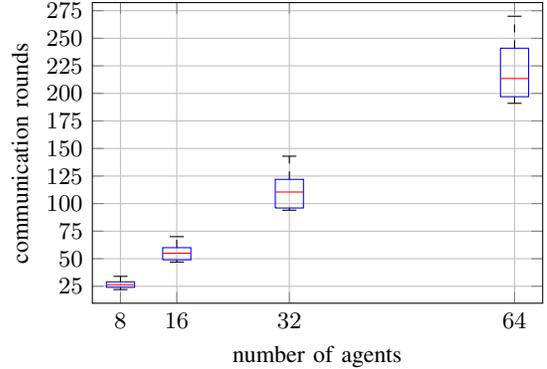

Next, we perform a numerical Monte Carlo analysis of the algorithm convergence while varying the network size and
its diameter. We recall that, in a static digraph, the diameter is the maximum distance taken all over the pair of
agents ($i,j$), where the distance is defined as the length of the shortest
directed path from $i$ to $j$. In the following we denote the diameter by $\diam$. 
We choose a cyclic digraph for which the diameter is proportional to the number of
agents, specifically $\diam=N_{ag}-1$. %
In particular, we consider the following cases: number of agents equal to $8$, $16$, $32$ and $64$. 
For each case, we generate $50$ random MILPs ensuring that each test case has a
non-empty set of feasible solutions.
The results are shown in Figure~\ref{fig:boxplot}. 
The red center line of each box shows the median value of the communication rounds for the $50$ random MILPs with fixed diameter. 
We highlight that the number of communication rounds needed for the convergence
grows linearly with the graph
diameter. %

\section{Conclusion}
\label{sec:concl}
In this paper, we proposed a distributed algorithm to solve Mixed Integer Linear
Programs over peer-to-peer networks. 
In the proposed distributed setup, the constraints of the \ac{MILP} are assigned to a
network of agents. The agents have a limited amount of memory and computation
capabilities and are able to communicate with neighboring agents.
Following the idea of centralized cutting plane methods for \ac{MILP}, each
agent solves local (LP) relaxations of the global problem, generates cutting
planes, and exchanges active constraints (a candidate basis) with neighbors.
We proved that agents reach consensus on the lex-optimal solution of the
\ac{MILP} in a finite number of communication rounds.
Finally, we performed a set of numerical computations suggesting that the completion time
of the algorithm scales nicely with the number of agents..

\appendix[Mixed Integer Gomory cutting planes]
\renewcommand{\thesection}{A}
\label{sec:MIGCP}
\numberwithin{equation}{section}
\numberwithin{theorem}{section}
In order to derive a \ac{MIG} cut, \cite{gomory1960algorithm}, for a generic LP relaxation
of problem \eqref{eq:MILP}, let 
$\var_{LP}$ be the current lex-optimal solution and $B_{LP}$ an associated basis. 
From the definition of basis, the lex-optimal solution can be obtained by solving the following LP problem
\begin{equation}\label{eq:basisLP}
\begin{split}
\min_{\var} &\; c^\top \var\\ 
\subj &\; A_B \var \leq b_B \\
\end{split}
\end{equation}
where $A_B\in\real^{d\times d}$ and $b_B\in \real^d$ are the matrices obtained
by writing in vector form the inequalities associated to the basis $B_{LP}$.
We proceed by rewriting problem~\eqref{eq:basisLP} in standard form. 
As described in \cite{luenberger1984linear}, we i)
reformulate~\eqref{eq:basisLP} as a maximization problem, ii) replace $z$ with
two new decision variables $z_+ \in \real^{d}$ and $z_- \in \real^{d}$ having
nonnegative components, such that $z=z_+ - z_-$, iii) introduce positive slack
variables $s \in \real^{d}$. In the new set of variables, we have
\begin{equation}\label{eq:SFMILP}
	\begin{aligned}
	 	\underset{\xsf}{\text{max}} \,\,  \bar{c}^\top \xsf \\
		\subj \,\, \bar{A} \xsf &= {b_B}\\
		\xsf_i &\geq 0, k = 1,\ldots,3d \\
	\end{aligned}
\end{equation}
where $\xsf = [z_+^\top, z_-^\top, s^\top]^\top$,
$\bar{A}= [A_B, -A_B, I_{d}]$,
$\bar{c}= [-c^\top, c^\top, 0_{d}^\top]^\top$, $I_{d}$ is the
identity matrix of dimension $d \times d$, $0_{d}$ is $d$-dimensional zero
vector.
Let the matrix $\bar{A}$ be partitioned as $\bar{A} = [\bar{A}_B, \bar{A}_N]$
where $\bar{A}_B$ is a suitable\footnote{If $\var_{LP_k}\geq 0$
  $\forall k=1\ldots d$, then $\bar{A}_B=A_B$. If $\var_{LP_k}\leq 0$
  $\forall k=1\ldots d$, then $\bar{A}_B=-A_B$. In the general case $\bar{A}_B$
  is a suitable permutation of columns of $A_B$ and $-A_B$.} $d \times d$
non-singular submatrix of $\bar{A}$ and $\bar{A}_N$ consists of the remaining
columns of $\bar{A}$.
Define the corresponding partition of the vector $u = [u_B^\top, u_N^\top]^\top$ (basic and nonbasic variables). 
The basic solution corresponding to the basis matrix $\bar{A}_B$ is given by
\begin{equation} \label{eq:SimplexTableau_red}
	\xsf_{B} = \bar{b} - \bar{a} u_N \,,
\end{equation}
where $\bar{a} = \bar{A}_B^{-1}\bar{A}_N$ and $\bar{b} = \bar{A}_B^{-1}b_B$. 
The \ac{MIG} cut is derived from the row of the simplex tableau \eqref{eq:SimplexTableau_red} corresponding to 
a basic variable, let us say the $k$-th component of~\eqref{eq:SimplexTableau_red}, that is required to be integer but it is not in the current solution. 
For the $k$-th row of the tableau, let us define $N_+:=\{\ell: \bar{a}_{k\ell} \geq 0\}$ and $N_-:=\{\ell:\bar{a}_{k\ell} < 0\}$,
where $\bar{a}_{k\ell}$ is the entry of the $(k,\ell)$-th entry of matrix $\bar{a}$. 
Then the \ac{MIG} cut is
\begin{align}
\label{eq:MILPCut}
\begin{split}
  &\sum_{\substack{f_{k\ell}\leq f_0\\\ell \text{
        integer}}}f_{k\ell}u_\ell+\frac{f_0}{1-f_0} \sum_{\substack{f_{k\ell}> f_0\\\ell
      \text{ integer}}}(1-f_{k\ell})u_\ell +\\ 
  &+\sum_{\substack{\ell\in N_+\\\ell \text{
        non-integer}}}{\bar{a}}_{k\ell}u_\ell-\frac{f_0}{1-f_0} \sum_{\substack{\ell\in
      N_-\\\ell \text{ non-integer}}}{\bar{a}}_{k\ell}u_\ell \geq f_0,\\ 
\end{split}
\end{align}
where $f_{k\ell}$ is the fractional part of $\bar{a}_{k\ell}$, and $f_0$ is the fractional part of the $k$-th component of $A_B^{-1}b_B$.
We can rewrite the \ac{MIG} cut \eqref{eq:MILPCut} with respect to the original
decision variables $\var$ (as in~\eqref{eq:MILPCut_z}) by taking in mind that
$s=b_B-A_B(\var_+-\var_-)$ and $\var=\var_+-\var_-$.

\begin{theorem}[\cite{jorg2008k}]
	The \ac{MIG} cut~\eqref{eq:MILPCut} is a valid cutting plane for $S =\{P \cap \integer^{d_Z} \times \real^{d_R}\}$. 
\end{theorem}

\begin{theorem}[\cite{balas2013combining}]%
\label{thm:MIG_intcut}
Let $u^\star$ be a solution of \eqref{eq:SFMILP} and $B$ the associated
basis. Let us suppose the $k$-th element of $u^\star$ is not
integer. %
Then the intersection cut to the split disjunction $D(e_k, \floor{u_k^\star} )$
and the basis $B$ is equal to the \ac{MIG} cut to the $k$-th component. 
\end{theorem}


\end{document}